\documentclass[12pt]{article}

\usepackage{latexsym,amssymb,amsmath,amsfonts,amsthm}
\usepackage{float}
\usepackage{xcolor}
\usepackage{geometry} 

\newtheorem{defi}{Definition}

\newtheorem{tm}{Theorem}

\newtheorem{rem}{Remark}

\begin{document}

\vspace*{0.5cm}

\begin{center}
{\Large\bf  Distance-regular graphs obtained from the Mathieu groups} 
\end{center}

\vspace*{0.5cm}

\begin{center}
    Dean Crnkovi\'c (deanc@math.uniri.hr)\\
		Nina Mostarac (nmavrovic@math.uniri.hr)\\
		and \\
		Andrea \v Svob (asvob@math.uniri.hr)\\[3pt]
		{\it\small Department of Mathematics} \\
		{\it\small University of Rijeka} \\
		{\it\small Radmile Matej\v ci\'c 2, 51000 Rijeka, Croatia}\\
\end{center}

\vspace*{0.5cm}

\begin{abstract}
In this paper we construct distance-regular graphs admitting a transitive action of the five sporadic simple groups discovered by E. Mathieu, 
the Mathieu groups $M_{11}$, $M_{12}$, $M_{22}$, $M_{23}$ and $M_{24}$. From the code spanned by the adjacency matrix of the strongly regular graph with parameters (176,70,18,34) we obtain block designs having the full automorphism groups isomorphic to the Higman-Sims finite simple group.
Further, we discuss a possibility of permutation decoding of the codes spanned by the adjacency matrices of the graphs constructed and find small PD-sets for some of the codes.
\end{abstract}

\bigskip

\noindent {\bf AMS classification numbers:} 05E18, 05E30, 94B05, 05B05.


\section{Introduction}

The main motivation for this paper is to give further contribution to the classification of transitive distance regular graphs (DRGs), especially those admitting a transitive action of a simple group. 
The research presented in the paper can be seen as a continuation of the work given in \cite{cms-m11}, in which transitive structures constructed from the Mathieu group $M_{11}$ were described.
On the other hand, this is also a continuation of the work by Praeger and Soicher \cite{low-rank}, where they study the graphs admitting a sporadic simple group or its automorphism group as a 
vertex-transitive group of automorphisms of rank at most 5. 
However, in this paper we are also interested in graphs admitting a vertex-transitive automorphism groups of ranks greater than 5.
In the literature, there are also other examples of the research in design and graph theory in which finite groups had important contribution (see \cite{crc-srg, drg, crsSRG, crsMANY, moori}).
In this paper we are focused on DRGs admitting a transitive action of one of the five Mathieu sporadic simple groups.

For relevant background reading in the group theory we refer the reader to \cite{atlas, r}, to \cite{crc-srg, tonchev-book} for the theory of strongly regular graphs, 
to \cite{drg,dam} for the theory of distance-regular graphs, and to \cite{Huffman, MacS} for coding theory. 

In this paper we study the Mathieu groups $M_{11}$, $M_{12}$, $M_{22}$, $M_{23}$ and $M_{24}$, which are the sporadic simple groups of orders 
7920, 95040, 443520, 10200960 and 244823040, respectively. The Mathieu groups $M_{11}$, $M_{12}$, $M_{22}$, $M_{23}$ and $M_{24}$ were found by Emile Mathieu in \cite{mathieu1, mathieu2, mathieu3}. 
They are multiply transitive permutation groups on 11, 12, 22, 23 or 24 objects, and they were the first sporadic groups discovered (see \cite{wielandt, wilson}).
The Mathieu groups have been studied so far in various literature and in various settings. For example,
$t$-designs arising from Mathieu groups $M_{22}$, $M_{23}$ and $M_{24}$ have been studied in a work by Kramer, Magliveras and Mesner (see \cite{kra-mag-mes}) and those arising from $M_{11}$ in 
\cite{cms-m11}. Codes connected with Mathieu group $M_{11}$ have been studied in \cite{vmc-novak} and those connected with $M_{12}$ in \cite{bailey}. More about small representations (up to rank 5) 
of these five sporadic simple groups can be found in \cite{low-rank}. We refer the reader to \cite{atlas, wilson} for more details about these groups. Using the method outlined in Section 
\ref{SRG_groups}, in this paper we constructed and partially classified DRGs (SRGs and DRGs of diameter $d \ge 3$) from above mentioned simple groups.


We also study codes spanned by the adjacency matrices of the constructed DRGs. 
Codes with large automorphism groups are suitable for permutation decoding (see \cite{key,key-pd, perm-mcw}), the decoding method developed by Jessie MacWilliams in the early 60s 
that can be used when a linear code has a sufficiently large automorphism group to ensure the existence of a set of automorphisms, called a PD-set, that has some specific properties.
Therefore, the codes constructed in this paper are suitable for permutation decoding, and we found small PD-sets for some of the constructed codes.
From the adjacency matrix of the SRG with parameters (176,70,18,34), beside the well-known symmetric block design with parameters (176,50,14) having the full automorphism group isomorphic to the Higman-Sims group, we obtain block designs with parameters 2-(176,56,110), 2-(176,64,540), 2-(176,66,780), 2-(176,70,2760), 2-(176,72,2556), 2-(176,78,37752), 2-(176,80,124030), 2-(176,82,99630), 
2-(176,86,87720) and 2-(176,88,210540) as support designs. The full automorphism groups of the new block designs are
also isomorphic to the Higman-Sims group HS. Further, one of the constructed SRGs with parameters (144,55,22,20) is new, up to our best knowledge.
 
To find the graphs and compute their full automorphism groups, and to obtain PD-sets and the corresponding information sets of the codes,  we used programmes written for Magma \cite{magma} and GAP \cite{GAP2020}. The constructed DRGs, including the SRGs, and the obtained PD-sets and the corresponding information sets of the codes
can be found at the link:

\begin{verbatim}
http://www.math.uniri.hr/~asvob/DRGs_Mathieu_PD.7z.
\end{verbatim} 

\section{Preliminaries}

In this section we begin with the required definitions and notation.

\begin{defi}\label{CohConf}
A coherent configuration on a finite non-empty set $\Omega$ is an ordered pair $(\Omega, \mathcal{R})$ with $\mathcal{R}=\{R_0,R_1,\dots,R_d\}$ a set of non-empty relations on $\Omega$, 
such that the following axioms hold.

\begin{itemize}
	\item[(i)] $\displaystyle\sum_{i=0}^t R_i$ is the identity relation, where $\{R_0,R_1,\dots,R_t\} \subseteq \{R_0,R_1,\dots,R_d\}$.
		\item[(ii)] $\mathcal{R}$ is a partition of $\Omega^2$.
		\item[(iii)] For every relation $R_i\in \mathcal{R}$, its converse $R_i^T=\{(y,x):(x,y)\in R_i\}$ is in $\mathcal{R}$.
		\item[(iv)] There are constants $p_{ij}^k$ known as the intersection numbers of the coherent configuration $\mathcal{R}$, such that for $(x,y)\in R_k$, the number of elements $z$ in $\Omega$ 
		for which $(x,z)\in R_i$ and $(z,y)\in R_j$ equals $p_{ij}^k$.
\end{itemize}
\end{defi}

We say that a coherent configuration is homogeneous if it contains the identity relation, i.e., if $R_0=I$.
If $\mathcal{R}$ is a set of symmetric relations on $\Omega$, then a coherent configuration is called symmetric. 
A symmetric coherent configuration is homogeneous (see \cite{cameron}).
Symmetric coherent configurations are introduced by Bose and Shimamoto in \cite{bose-associate} and called association schemes.
An association scheme with relations $\{R_0,R_1,\dots,R_d\}$ is called a $d$-class association scheme.

Let $\Gamma$ be a graph with diameter $d$, and let $\delta(u,v)$ denote the distance between vertices $u$ and $v$ of $\Gamma$.
The $i$th-neighborhood of a vertex $v$ is the set $\Gamma_{i}(v) = \{ w : \delta(v,w) = i \}$. 
Similarly, we define $\Gamma_{i}$ to be the $i$th-distance graph of $\Gamma$, that is, the vertex set of $\Gamma_{i}$ is the same as for $\Gamma$, with adjacency in $\Gamma_{i}$ defined by the $i$th distance relation in $\Gamma$.
We say that $\Gamma$ is distance-regular if the distance relations of $\Gamma$ give the relations of a $d$-class association scheme, that is, for every choice of $0 \leq i,j,k \leq d$, all vertices $v$ and $w$ with $\delta(v,w)=k$ satisfy $|\Gamma_{i}(v) \cap \Gamma_{j}(w)| = p^{k}_{ij}$ for some constant $p^{k}_{ij}$.
In a distance-regular graph, we have that $p^{k}_{ij}=0$ whenever $i+j < k$ or $k<|i-j|$.
A distance-regular graph $\Gamma$ is necessarily regular with degree $p^{0}_{11}$; more generally, each distance graph $\Gamma_{i}$ is regular with degree $k_{i}=p^{0}_{ii}$.

An equivalent definition of distance-regular graphs is the existence of the constants $b_{i}=p^{i}_{i+1,1}$ and $c_{i}= p^{i}_{i-1,1}$ for $0 \leq i \leq d$ (notice that $b_{d}=c_{0}=0$).
The sequence $\{b_0,b_1,\dots,b_{d-1};c_1,c_2,\dots,c_d\}$, where $d$ is the diameter of $\Gamma$ is called the intersection array of $\Gamma$. 
Clearly, $b_0=k$, $b_d=c_0=0$, $c_1=0$.

A regular graph is strongly regular with parameters $(v,k, \lambda , \mu )$ if it has $v$ vertices, degree $k$,
and if any two adjacent vertices are together adjacent to $\lambda$ vertices, while any two non-adjacent vertices are together adjacent to $\mu$ vertices.
A strongly regular graph with parameters $(v,k, \lambda , \mu )$ is usually denoted by SRG$(v, k, \lambda, \mu)$. 
A strongly regular graph is a distance-regular graph with diameter 2 whenever $\mu \neq 0$.
The intersection array of an SRG is given by $\{k,k-1-\lambda;1,\mu\}$.

\section{DRGs constructed from the Mathieu groups} \label{SRG_groups}

Let $G$ be a finite permutation group acting on the finite set $\Omega$. This action induce the action of the group $G$ on the set $\Omega\times \Omega$. For more information see \cite{wielandt}. 
The orbits of this action are the sets of the form $\{(\alpha g,\beta g): g\in G\}$. If $G$ is transitive, then $\{(\alpha ,\alpha):\alpha\in \Omega\}$ is one such orbit. 
If the rank of $G$ is $r$, then it has $r$ orbits on $\Omega\times\Omega$.
Let $|\Omega|=n$ and $\Delta_i$ is one of these orbits. We say that the $n\times n$ matrix $A_i$, with rows and columns indexed by $\Omega$ and entries 

$$
A_i(\alpha,\beta)= \left \lbrace
\begin{array}{ll}
1,      & \mathrm{if} \ (\alpha,\beta)\in \Delta_i\\
0, & \mathrm{otherwise}.
\end{array}  \right.
$$
is called the adjacency matrix for the orbit $\Delta_i$.

Let $A_0,\dots,A_{r-1}$ be the adjacency matrices for the orbits of $G$ on $\Omega\times\Omega$. These satisfy the following conditions.

\begin{itemize}
\item[(i)] $A_0=I$, if $G$ is transitive on $\Omega$. If $G$ has $s$ orbits on $\Omega$, then $I$ is a sum of $s$ adjacency matrices.
\item[(ii)] $\displaystyle\sum_{i}A_i=J$, where $J$ is the all-one matrix.
\item[(iii)] If $A_i$ is an adjacency matrix, then so is its transpose $A_i^T$.
\item[(iv)] If $A_i$ and $A_j$ are adjacency matrices, then their product is an integer-linear-combination of adjacency matrices.
\end{itemize}

If $A_i$ is symmetric, then the corresponding orbit is called self-paired. Further, if $A_i=A_j^T$, then the corresponding orbits are called mutually paired.

The graphs obtained in this paper are constructed using the method described in \cite{cms} which can be rewritten in terms of coherent configurations in the following way.

\begin{tm}\label{main}
Let $G$ be a finite permutation group acting transitively on the set $\Omega$ and $A_0,\dots,A_d$ be the adjacency matrices for orbits of $G$ on $\Omega\times\Omega$. 
Let $\{ B_1,\dots, B_t \} \subseteq \{ A_1,\dots, A_d \}$ be a set of adjacency matrices for the self-paired or mutually paired orbits. 
Then $\displaystyle M=\sum_{i=1}^{t}B_i$ is the adjacency matrix of a regular graph $\Gamma$. The group $G$ acts transitively on the set of vertices of the graph $\Gamma$.
\end{tm}

Using this method one can construct all regular graphs admitting a transitive action of the group $G$. We will be interested only in those 
regular graphs that are distance-regular, and especially strongly regular. 

\begin{rem}\label{running_time}
Because of the large number of possibilities for building the first row of the adjacency matrix of a DRG, the only way to obtain the classification of DRGs given in this paper was with the use of computers. The running time complexity of the algorithm used for the construction of graphs depends on a number of parameters, such as the size of the used subgroup, the number of orbits of a vertex stabilizer, the number of vertices of the graphs and the number of self-paired and mutually paired orbits in a particular case. 
\end{rem}

\subsection{DRGs from the group $M_{11}$}\label{M11}

The Mathieu group $M_{11}$ has the order 7920 and up to conjugation has 39 subgroups. 
In Table \ref{tb:subgrpsM11} we give the list of all the subgroups $H_i^1 \le M_{11}$ which lead to the construction of SRGs or DRGs of diameter $d \ge 3$.

\begin{table}[H]
\begin{center} \begin{scriptsize}
\begin{tabular}{|c|c|r|r|r|c|}
\hline
Subgroup& Structure& Order & Index & Rank & Primitive\\
\hline
\hline
$H^1_{1}$ & $S_5$ & 144 & 55 & 3 &  yes\\
$H^1_{2}$ & $Z_9:QD16$ & 120 & 66 & 4 & yes \\
$H^1_{3}$ & $Z_{11}:Z_5$ & 55 & 144 & 6 & no \\
$H^1_{4}$ & $GL(2,3)$ &  48& 165 & 8 &  yes\\
$H^1_{5}$ & $S_3\times S_3$ & 36 & 220 & 16 & no \\
$H^1_{6}$ & $S_4$ & 24 & 330 & 23 & no \\

\hline
\hline
\end{tabular}\end{scriptsize} 
\caption{\footnotesize Subgroups of the group $M_{11}$}\label{tb:subgrpsM11}
\end{center}
\end{table}

Using the method described in Theorem \ref{main} we obtained all DRGs with at most 2000 vertices and for which the rank of the permutation representation of the group is at most 25, i.e. we gave the classification of such DRGs.

\begin{tm} \label{srg-M11}
Up to isomorphism there are exactly five strongly regular graphs and exactly three distance-regular graphs of diameter $d \ge 3$ with at most 2000 vertices and for which the rank of the permutation representation of the group is at most 25, 
admitting a transitive action of the group $M_{11}$.
The SRGs have parameters $(55,18,9,4)$, $(66,20,10,4)$, $(144,55,22,20)$, $(144,66,30,30)$ and $(330,63,24,9)$, 
and the DRGs have $165$, $220$ and $330$ vertices, respectively.
Details about the obtained strongly regular graphs are given in Table \ref{tb:srgM11} and details about the obtained DRGs with $d \ge 3$ are given in Table \ref{tb:drgsM11}.
\end{tm}

\begin{table}[H]
\begin{center} \begin{scriptsize}
\begin{tabular}{|c|c|c|}
\hline
Graph $\Gamma$ & Parameters  & $Aut (\Gamma) $  \\
\hline
\hline
$\Gamma^{1}_{1}=\Gamma(M_{11},H^1_{1})$ & (55,18,9,4) & $S_{11}$  \\ 
$\Gamma^{1}_{2}=\Gamma(M_{11},H^1_{2})$ & (66,20,10,4) & $S_{12}$  \\ 
$\Gamma^{1}_{3}=\Gamma(M_{11},H^1_{3})$ & (144,55,22,20) & $M_{11}$  \\ 
$\Gamma^{1}_{4}=\Gamma(M_{11},H^1_{3})$ & (144,66,30,30) & $M_{12}:Z_2$  \\ 
$\Gamma^{1}_{5}=\Gamma(M_{11},H^1_{6})$ & (330,63,24,9) & $S_{11}$  \\ 
\hline
\hline
\end{tabular}\end{scriptsize} 
\caption{\footnotesize SRGs constructed from the group $M_{11}$}\label{tb:srgM11}
\end{center}
\end{table}

\begin{table}[H]
\begin{center} \begin{scriptsize}
\begin{tabular}{|c|c|c|c|c|}
\hline
Graph $\Gamma$ & Number of vertices & Diameter  & Intersection array & $Aut (\Gamma) $  \\
\hline
\hline
$\Gamma^{1}_{6}=\Gamma(M_{11},H^1_{4})$ &  165 & 3 & $\{ 24,14,6;1,4,9 \}$& $S_{11}$\\ 
$\Gamma^{1}_{7}=\Gamma(M_{11},H^1_{5})$ &  220 & 3 & $\{ 27,16,7;1,4,9 \}$& $S_{12}$\\ 
$\Gamma^{1}_{8}=\Gamma(M_{11},H^1_{6})$ &  330 & 4 & $\{ 28,18,10,4;1,4,9,16 \}$& $S_{11}$\\ 
\hline
\hline
\end{tabular} \end{scriptsize}
\caption{\footnotesize DRGs constructed from the group $M_{11}$, $d \geq$ 3}\label{tb:drgsM11}
\end{center}
\end{table}

{\bf Proof}.
There are 39 conjugacy classes of subgroups of $M_{11}$, but only 19 of them lead to a permutation representation of rank at most 25 and of index at most 2000. Applying the method described in Theorem \ref{main} to the permutation representations on cosets of these 19 subgroups we obtain the results.
{$\Box$}

\begin{rem} \label{IsoSrgsM11} 
All SRGs given in Table \ref{tb:srgM11} are isomorphic to the ones constructed in \cite{cms-m11}.
Since the SRG $\Gamma_3^1$ does not yield a partial geometry, it cannot be obtained from an orthogonal array.
\end{rem}

\begin{rem}
The graphs $\Gamma_6^1$, $\Gamma_7^1$ and $\Gamma_8^1$ are unique graphs with the given intersection arrays, known as Johnson graphs, $J(11,3)$, $J(12,3)$ and $J(11,4)$, respectively (see \cite{drg}).
\end{rem}

\subsection{DRGs from the group $M_{12}$}\label{M12}

The Mathieu group $M_{12}$ has the order 95040 and up to conjugation has 147 subgroups. 
In Table \ref{tb:subgrpsM12} we give the list of all the subgroups $H_i^2 \le M_{12}$ which lead to the construction of SRGs or DRGs of diameter $d \ge 3$.

\begin{table}[H]
\begin{center} \begin{scriptsize}
\begin{tabular}{|c|c|r|r|r|c|}
\hline
Subgroup& Structure& Order & Index & Rank & Primitive\\
\hline
\hline
$H^2_{1}$ & $(A_6.Z_2):Z_2$ & 1440 & 66 & 3 &  yes\\
$H^2_{2}$ & $L(2,11)$ & 660 & 144 & 5 & no \\
$H^2_{3}$ & $((E_9:Q_8):Z_3):Z_2$ & 432 & 220 & 5 & yes \\
$H^2_{4}$ & $((E_8:E_4):Z_3):Z_2$ &  192& 495 & 11 &  yes\\
$H^2_{5}$ & $S_5$ & 120 & 792 & 15 & no \\
\hline
\hline
\end{tabular}\end{scriptsize} 
\caption{\footnotesize Subgroups of the group $M_{12}$}\label{tb:subgrpsM12}
\end{center}
\end{table}

Using the method described in Theorem \ref{main} we obtained all DRGs with at most 2000 vertices and for which the rank of the permutation representation of the group is at most 20, i.e. we gave the classification of such DRGs.

\begin{tm} \label{srg-M12}
Up to isomorphism there are exactly seven strongly regular graphs and exactly three distance-regular graphs of diameter $d \ge 3$ with at most 2000 vertices and for which the rank of the permutation representation of the group is at most 20, 
admitting a transitive action of the group $M_{12}$.
The SRGs have parameters $(66,20,10,4)$, $(144,66,30,30)$, $(144,55,22,20)$, $(144,22,10,2)$ and $(495,238,109,119)$, 
and the DRGs have $220$, $495$ and $792$ vertices, respectively.
Details about the obtained strongly regular graphs are given in Table \ref{tb:srgM12} and details about the obtained DRGs with $d \ge 3$ are given in Table \ref{tb:drgsM12}.
\end{tm}

\begin{table}[H]
\begin{center} \begin{scriptsize}
\begin{tabular}{|c|c|c|}
\hline
Graph $\Gamma$ & Parameters  & $Aut (\Gamma) $  \\
\hline
\hline
$\Gamma^{2}_{1}=\Gamma(M_{12},H^2_{1})$ & (66,20,10,4) & $S_{12}$  \\ 
$\Gamma^{2}_{2}=\Gamma(M_{12},H^2_{2})$ & (144,66,30,30) & $M_{12}:Z_2$  \\ 
$\Gamma^{2}_{3}=\Gamma(M_{12},H^2_{2})$ & (144,66,30,30) & $M_{12}:Z_2$  \\ 
$\Gamma^{2}_{4}=\Gamma(M_{12},H^2_{2})$ & (144,66,30,30) & $M_{12}$  \\ 
$\Gamma^{2}_{5}=\Gamma(M_{12},H^2_{2})$ & (144,55,22,20) & $M_{12}:Z_2$  \\ 
$\Gamma^{2}_{6}=\Gamma(M_{12},H^2_{2})$ & (144,22,10,2) & $S_{12}\wr S_2$  \\ 
$\Gamma^{2}_{7}=\Gamma(M_{12},H^2_{4})$ & (495,238,109,119) & $O^{-}(10,2):Z_2$  \\ 
\hline
\hline
\end{tabular}\end{scriptsize} 
\caption{\footnotesize SRGs constructed from the group $M_{12}$}\label{tb:srgM12}
\end{center}
\end{table}

\begin{table}[H]
\begin{center} \begin{scriptsize}
\begin{tabular}{|c|c|c|c|c|}
\hline
Graph $\Gamma$ & Number of vertices & Diameter  & Intersection array & $Aut (\Gamma) $  \\
\hline
\hline
$\Gamma^{2}_{8}=\Gamma(M_{12},H^2_{3})$ &  220 & 3 & $\{ 27,16,7;1,4,9 \}$& $S_{12}$\\ 
$\Gamma^{2}_{9}=\Gamma(M_{12},H^2_{4})$ &  495 & 4 & $\{ 32,21,12,5;1,4,9,16 \}$& $S_{12}$\\ 
$\Gamma^{2}_{10}=\Gamma(M_{12},H^2_{5})$ &  792 & 5 & $\{ 35,24,15,8,3;1,4,9,16,25 \}$& $S_{12}$\\ 
\hline
\hline
\end{tabular} \end{scriptsize}
\caption{\footnotesize DRGs constructed from the group $M_{12}$, $d \geq$ 3}\label{tb:drgsM12}
\end{center}
\end{table}

{\bf Proof}.
There are 147 conjugacy classes of subgroups of $M_{12}$, but only 31 of them lead to a permutation representation of rank at most 20 and of index at most 2000. Applying the method described in Theorem \ref{main} to the permutation representations on cosets of these 31 subgroups we obtain the results.
{$\Box$}

\begin{rem} \label{IsoSrgsM12}
The strongly regular graph $\Gamma_1^2$ is isomorphic to the triangular graph $T(12)$. The adjacency matrices of non-isomorphic SRGs $\Gamma^{2}_{2}$, $\Gamma^{2}_{3}$ and $\Gamma^{2}_{4}$ are the incidence matrices of symmetric designs with parameters $(144,66,30)$, designs with Menon parameters (related to a regular Hadamard matrix of order 144).
These symmetric designs have been described in \cite{lempken, wirth}.  According to Brouwer's table (see  \cite{aeb}), known graphs with the parameters equal to the parameters of the graph $\Gamma_5^2$ (not isomorphic to $\Gamma_3^1$) are obtainable from orthogonal arrays $OA(12,5)$. Since the SRG $\Gamma_5^2$ cannot be obtained from orthogonal array, our graph is new. So far, according to \cite{aeb}, the graphs $\Gamma_3^1$ and $\Gamma_5^2$ are the only known graphs with these parameters not arising from orthogonal array. The graph $\Gamma_{6}^2$ is unique graph with the given parameters and the graph $\Gamma_{7}^2$ is isomorphic to the $O^{-}(10,2)$ polar graph.
Strongly regular graphs with parameters $(144,66,30,30)$ have been known before (see \cite{crc-srg,aeb}). 

\end{rem}

\begin{rem}
The graphs $\Gamma_8^2$, $\Gamma_9^2$ and $\Gamma_{10}^2$ are unique graphs with the given intersection arrays, known as Johnson graphs, $J(12,3)$, $J(12,4)$ and $J(12,5)$, respectively (see \cite{drg}).
\end{rem}

\subsection{DRGs from the group $M_{22}$}\label{SRG_M22}

The Mathieu group $M_{22}$ has the order 443520 and up to conjugation 156 subgroups. 
In Table \ref{tb:subgrpsM22} we give the list of all the subgroups $H_i^3 \le M_{22}$ which lead to the construction of SRGs or DRGs of diameter $d \ge 3$.

\begin{table}[H]
\begin{center} \begin{scriptsize}
\begin{tabular}{|c|c|r|r|r|c|}
\hline
Subgroup& Structure& Order & Index & Rank & Primitive\\
\hline
\hline
$H^3_{1}$ & $E_{16}:A_6$& 5760 & 77 & 3 &  yes\\
$H^3_{2}$ & $A_7$& 2520 & 176 & 3 &  yes\\
$H^3_{3}$ & $E_{16}:S_5$& 1920 & 231 & 4 & yes \\
$H^3_{4}$ & $E_8:L(3,2)$ & 1344 & 330 & 5 & yes \\
$H^3_{5}$ & $L(2,11)$& 660  & 672 & 6 & yes \\
$H^3_{6}$ & $(A_4\times A_4):Z_2$ & 288 & 1540 & 22 &  no\\
\hline
\hline
\end{tabular}\end{scriptsize} 
\caption{\footnotesize Subgroups of the group $M_{22}$}\label{tb:subgrpsM22}
\end{center}
\end{table}

Using the method described in Theorem \ref{main} we obtained all DRGs with at most 2000 vertices and for which the rank of the permutation representation of the group is at most 30, i.e. we gave the classification of such DRGs.

\begin{tm} \label{srg-M22}
Up to isomorphism there are exactly five strongly regular graphs and exactly three distance-regular graphs of diameter $d \ge 3$ with at most 2000 vertices and for which the rank of the permutation representation of the group is at most 30, 
admitting a transitive action of the group $M_{22}$.
The SRGs have parameters $(77,16,0,4)$, $(176,70,18,34)$, $(231,30,9,3)$, $(231,40,20,4)$ and $(672,176,40,48)$, 
and the DRGs have $330$, $672$ and $1540$ vertices, respectively.
Details about the obtained strongly regular graphs are given in Table \ref{tb:srgM22} and details about the obtained DRGs with $d \ge 3$ are given in Table \ref{tb:drgsM22}.
\end{tm}

\begin{table}[H]
\begin{center} \begin{scriptsize}
\begin{tabular}{|c|c|c|}
\hline
Graph $\Gamma$ & Parameters  & $Aut (\Gamma) $  \\
\hline
\hline
$\Gamma^{3}_{1}=\Gamma(M_{22},H^3_{1})$ & $(77,16,0,4)$ & $M_{22}:Z_2$  \\ 
$\Gamma^{3}_{2}=\Gamma(M_{22},H^3_{2})$ & $(176,70,18,34)$ & $M_{22}$  \\ 
$\Gamma^{3}_{3}=\Gamma(M_{22},H^3_{3})$ & $(231,30,9,3)$ & $M_{22}:Z_2$  \\ 
$\Gamma^{3}_{4}=\Gamma(M_{22},H^3_{3})$ & $(231,40,20,4)$ & $S_{22}$  \\ 
$\Gamma^{3}_{5}=\Gamma(M_{22},H^3_{5})$ & $(672,176,40,48)$ & $(U(6,2):Z_2):Z_2$  \\ 
\hline
\hline
\end{tabular}\end{scriptsize} 
\caption{\footnotesize SRGs constructed from the group $M_{22}$}\label{tb:srgM22}
\end{center}
\end{table}

\begin{table}[H]
\begin{center} \begin{scriptsize}
\begin{tabular}{|c|c|c|c|c|}
\hline
Graph $\Gamma$ & Number of vertices& Diameter  & Intersection array & $Aut (\Gamma) $  \\
\hline
\hline
$\Gamma^{3}_{6}=\Gamma(M_{22},H^3_{4})$ &  330 & 4 & $\{ 7,6,4,4; 1,1,1,6 \}$& $M_{22}:Z_2$\\ 
$\Gamma^{3}_{7}=\Gamma(M_{22},H^3_{5})$ &  672 & 3 & $\{ 110,81,12; 1,18,90 \}$& $M_{22}:Z_2$\\ 
$\Gamma^{3}_{8}=\Gamma(M_{22},H^3_{6})$ &  1540 & 3 & $\{ 57,36,17; 1,4,9 \}$& $S_{22}$\\ 
\hline
\hline
\end{tabular} \end{scriptsize}
\caption{\footnotesize DRGs constructed from the group $M_{22}$, $d \geq$ 3}\label{tb:drgsM22}
\end{center}
\end{table}

{\bf Proof}.
There are 156 conjugacy classes of subgroups of $M_{22}$, but only 21 of them lead to a permutation representation of rank at most 30 and of index at most 2000. Applying the method described in Theorem \ref{main} to the permutation representations on cosets of these 21 subgroups we obtain the results.
{$\Box$}

\begin{rem} \label{IsoSrgsM22}
The strongly regular graphs $\Gamma^{3}_{1}$ and $\Gamma^{3}_{2}$ are unique graphs with these parameters. The graph $\Gamma^{3}_{3}$ is isomorphic to the SRG known as the Cameron graph. The SRG $\Gamma^{3}_{4}$ is isomorphic to the triangular graph $T(22)$ and $\Gamma^{3}_{5}$ is isomorphic to the graph known as $U(6,2)$-graph. For more information we refer the reader to \cite{crc-srg, aeb}.
\end{rem}

\begin{rem}
The graph $\Gamma^{3}_{6}$ is isomorphic to the graph known as $M_{22}$-graph or doubly truncated Witt graph. The graph $\Gamma^{3}_{7}$ is isomorphic to the one constructed by Soicher in \cite{soicher-drg}. So far, it is the only known example of DRG with this intersection array. The graph $\Gamma^{3}_{8}$ is known as Johnson graph $J(22,3)$.
(see \cite{drg})
\end{rem}

\subsection{DRGs from the group $M_{23}$}\label{SrgsM23}

The Mathieu group $M_{23}$ has order 10200960 and up to conjugation 204 subgroups. 
In Table \ref{tb:subgrpsM23} we give the list of all the subgroups $H_i^4 \le M_{23}$ which lead to the construction of SRGs or DRGs of diameter $d \ge 3$.

\begin{table}[H]
\begin{center} \begin{scriptsize}
\begin{tabular}{|c|c|r|r|r|c|}
\hline
Subgroup& Structure& Order & Index & Rank & Primitive\\
\hline
\hline
$H^4_{1}$ & $L(3,4):Z_2$ & 40320 & 253 & 3 & yes \\
$H^4_{2}$ & $E_{16}:A_7$&  40320& 253 &  3&  yes\\
$H^4_{3}$ & $A_8$& 20160 & 506 & 4 &  yes\\
$H^4_{4}$ & $M_{11}$& 7920& 1288 & 4 & yes \\
$H^4_{5}$ & $E_{16}:(A_5:S_3)$& 5760 & 1771 & 8 &  yes\\
\hline
\hline
\end{tabular}\end{scriptsize} 
\caption{\footnotesize Subgroups of the group $M_{23}$}\label{tb:subgrpsM23}
\end{center}
\end{table}

Using the method described in Theorem \ref{main} we obtained all DRGs with at most 10000 vertices and for which the rank of the permutation representation of the group is at most 20, i.e. we gave the classification of such DRGs.

\begin{tm} \label{srg-M23}
Up to isomorphism there are exactly three strongly regular graphs and exactly two distance-regular graphs of diameter $d \ge 3$ with at most 10000 vertices and for which the rank of the permutation representation of the group is at most 20, 
admitting a transitive action of the group $M_{23}$.
The SRGs have parameters $(253,42,21,4)$, $(253,112,36,60)$ and $(1288,495,206,180)$, 
and the DRGs have $506$ and $1771$ vertices, respectively.
Details about the obtained strongly regular graphs are given in Table \ref{tb:srgM23} and details about the obtained DRGs with $d \ge 3$ are given in Table \ref{tb:drgsM23}.
\end{tm}

\begin{table}[H]
\begin{center} \begin{scriptsize}
\begin{tabular}{|c|c|c|}
\hline
Graph $\Gamma$ & Parameters  & $Aut (\Gamma) $  \\
\hline
\hline
$\Gamma^{4}_{1}=\Gamma(M_{23},H^4_{1})$ & (253,42,21,4) & $S_{23}$  \\ 
$\Gamma^{4}_{2}=\Gamma(M_{23},H^4_{2})$ & (253,112,36,60) & $M_{23}$  \\ 
$\Gamma^{4}_{3}=\Gamma(M_{23},H^4_{4})$ & (1288,495,206,180) & $M_{24}$  \\ 
\hline
\hline
\end{tabular}\end{scriptsize} 
\caption{\footnotesize SRGs constructed from the group $M_{23}$}\label{tb:srgM23}
\end{center}
\end{table}

\begin{table}[H]
\begin{center} \begin{scriptsize}
\begin{tabular}{|c|c|c|c|c|}
\hline
Graph $\Gamma$ & Number of vertices  & Diameter & Intersection array & $Aut (\Gamma) $  \\
\hline
\hline
$\Gamma^{4}_{4}=\Gamma(M_{23},H^4_{3})$ & 506 & 3 & $\{ 15,14,12; 1,1,9 \}$& $M_{23}$\\ 
$\Gamma^{4}_{5}=\Gamma(M_{23},H^4_{5})$ & 1771 & 3 & $\{ 60,38,18; 1,4,9 \}$& $S_{23}$\\ 
\hline
\hline
\end{tabular} \end{scriptsize}
\caption{\footnotesize DRG constructed from the group $M_{23}$, $d \geq$ 3}\label{tb:drgsM23}
\end{center}
\end{table}

{\bf Proof}.
There are 204 conjugacy classes of subgroups of $M_{23}$, but only 14 of them lead to a permutation representation of rank at most 20 and of index at most 10000. Applying the method described in Theorem \ref{main} to the permutation representations on cosets of these 14 subgroups we obtain the results.
{$\Box$}

\begin{rem} \label{IsoSrgsM23}
The graph $\Gamma_{1}^4$ is isomorphic to the triangular graph $T(23)$. The graph $\Gamma^{4}_{2}$ can be constructed from the group $M_{23}$ as a rank $3$ graph, and $\Gamma^{4}_{3}$ (isomorphic to the graph $\Gamma_{2}^5$) can be constructed from the group $M_{24}$ as a rank $3$ graph.
\end{rem}

\begin{rem}
The graph $\Gamma^{4}_{4}$ is isomorphic to the distance-regular graph that can be obtained from residual design of Steiner system $S(5,8,24)$. The graph $\Gamma^{4}_{5}$ is known as Johnson graph $J(23,3)$. For more information we refer the reader to \cite{drg}.
\end{rem}

\subsection{DRGs from the Mathieu group $M_{24}$}\label{SRG_M24}

The Mathieu group $M_{24}$ has order 244823040  and up to conjugation 1529 subgroups. 
In Table \ref{tb:subgrpsM24} we give the list of all the subgroups $H_i^5 \le M_{24}$ which lead to the construction of SRGs or DRGs of diameter $d \ge 3$.

\begin{table}[H]
\begin{center} \begin{scriptsize}
\begin{tabular}{|c|c|r|r|r|c|}
\hline
Subgroup& Structure& Order & Index & Rank & Primitive\\
\hline
\hline
$H^5_{1}$ & $M_{22}:Z_2$&  887040& 276 & 3 & yes \\
$H^5_{2}$ & $E_{16}:A_8$& 322560 & 759 & 4 & yes \\
$H^5_{3}$ & $M_{12}:Z_2$& 190080 & 1288 & 3 &  yes\\
$H^5_{4}$ & $(L(3,4):Z_3):Z_2$& 120960 & 2024 & 5 & yes \\
\hline
\hline
\end{tabular}\end{scriptsize} 
\caption{\footnotesize Subgroups of the group $M_{24}$}\label{tb:subgrpsM24}
\end{center}
\end{table}

Using the method described in Theorem \ref{main} we obtained all DRGs with at most 10000 vertices and for which the rank of the permutation representation of the group is at most 20, i.e. we gave the classification of such DRGs.

\begin{tm} \label{srg-M24}
Up to isomorphism there are exactly two strongly regular graphs and exactly two distance-regular graphs of diameter $d \ge 3$ with at most 10000 vertices and for which the rank of the permutation representation of the group is at most 20, 
admitting a transitive action of the group $M_{24}$.
The SRGs have parameters $(276,44,22,4)$ and $(1288,495,206,180)$, 
and the DRGs have $759$ and $2024$ vertices, respectively.
Details about the obtained strongly regular graphs are given in Table \ref{tb:srgM24} and details about the obtained DRGs with $d \ge 3$ are given in Table \ref{tb:drgsM24}.
\end{tm}

\begin{table}[H]
\begin{center} \begin{scriptsize}
\begin{tabular}{|c|c|c|}
\hline
Graph $\Gamma$ & Parameters  & $Aut (\Gamma) $  \\
\hline
\hline
$\Gamma^{5}_{1}=\Gamma(M_{24},H^5_{1})$ & $(276,44,22,4)$ & $S_{24}$  \\ 
$\Gamma^{5}_{2}=\Gamma(M_{24},H^5_{3})$ & $(1288,495,206,180)$ & $M_{24}$  \\ 
\hline
\hline
\end{tabular}\end{scriptsize} 
\caption{\footnotesize SRGs constructed from the group $M_{24}$}\label{tb:srgM24}
\end{center}
\end{table}

\begin{table}[H]
\begin{center} \begin{scriptsize}
\begin{tabular}{|c|c|c|c|c|}
\hline
Graph $\Gamma$ & Number of vertices  & Diameter & Intersection array & $Aut (\Gamma) $  \\
\hline
\hline
$\Gamma^{5}_{3}=\Gamma(M_{24},H^5_{2})$ & 759 & 3 & $\{ 30,28,24; 1,3,15 \}$& $M_{24}$\\ 
$\Gamma^{5}_{4}=\Gamma(M_{24},H^5_{4})$ & 2024 & 3 & $\{ 63,40,19; 1,4,9 \}$& $S_{24}$\\ 
\hline
\hline
\end{tabular} \end{scriptsize}
\caption{\footnotesize DRGs constructed from the group $M_{24}$, $d \geq$ 3}\label{tb:drgsM24}
\end{center}
\end{table}

{\bf Proof}.
There are 1529 conjugacy classes of subgroups of $M_{24}$, but only 15 of them lead to a permutation representation of rank at most 20 and of index at most 10000. Applying the method described in Theorem \ref{main} to the permutation representations on cosets of these 15 subgroups we obtain the results.
{$\Box$}

\begin{rem} \label{IsoSrgsM24}
The graph $\Gamma_{1}^5$ is isomorphic to the triangular graph $T(24)$. The graph $\Gamma^{5}_{2}$ (isomorphic to the graph $\Gamma_{3}^4$) can be constructed from the group $M_{24}$ as a rank $3$ graph.
\end{rem}

\begin{rem}
The graph $\Gamma^{5}_{3}$ is unique distance-regular graph known as near hexagon which can be obtained from Steiner system $S(5,8,24)$. The graph $\Gamma^{5}_{4}$ is known as Johnson graph $J(24,3)$. For more information we refer the reader to \cite{drg}.
\end{rem}

\section{Codes}

A code $C$ of length $n$ over the alphabet $Q$ is a subset $C\subseteq Q^n$. Elements of the code are called codewords. A code $C$ is called a $p$-ary linear code of dimension 
$m$ if $Q = \mathbb{F}_p$, for a prime power $p$, and $C$ is an $m$-dimensional subspace of the vector space $(\mathbb{F}_p)^n$. For $Q = \mathbb{F}_2$ a code is called binary.

Let $x = (x_1,...,x_n)$ and $y = (y_1,...,y_n) \in \mathbb{F}_p^n$. The Hamming distance between the words $x$ and $y$ is the number $d(x,y) = |\{i : x_i \neq y_i\}|$. 
The minimum distance of the code $C$ is defined by $d = \min \{ d(x,y) : x, y\in C, x \neq y \}$. The weight of a codeword $x$ is $w(x) = d(x,0) = |\{i : x_i \neq 0\}|$. 
For a linear code the minimum distance equals the minimum weight $d = \min \{ w(x) : x \in C, x\neq 0 \}$.

A $p$-ary linear code of length $n$, dimension $k$, and minimum distance $d$ is called an $[n,k,d]_p$ code or $[n,k,d]$ code when the size $p$ of the field is not mentioned. A linear $[n,k,d]$ code can detect at most $d-1$ errors in one codeword and correct at most 
$t = \left\lfloor \frac{d-1}{2} \right\rfloor$ errors.
Two binary linear codes are isomorphic if one can be obtained from the other by permuting the coordinate positions. An automorphism of the code $C$ is an isomorphism from $C$ to $C$.  
The generator matrix of a code is a $k \times n$ matrix whose rows are the vectors of a base of the code. 
Every code is isomorphic to a code with the generator matrix in the standard form, i.e. in the form $[I_k,A]$, where $I_k$ is the identity
matrix of order $k$ and $A$ some $k\times (n-k)$ matrix.

\subsection{Permutation decoding}

Let $C\subseteq \mathbb{F}_p^n$ be a  linear $[n, k, d]$ code. For $I\subseteq \{1,...,n\}$, let $p_I:\mathbb{F}_p^n \to \mathbb{F}_p^{|I|}$, $x \mapsto x|_I$,
be the $I$-projection of $\mathbb{F}_p^n$. Then $I$ is called an information set for $C$ if $|I|=k$ and $p_I(C)=\mathbb{F}_p^{|I|}$.
The set of the first $k$ coordinates for a code with a generator matrix  in the standard form is an information set.
The first $k$ coordinates are then called information symbols and the last $n-k$ coordinates are the check symbols and they form the corresponding check set.

Let $C\subseteq \mathbb{F}_p^n$ be a linear $[n, k, d]$ code that can correct at most $t$ errors (i.e. $t$-error-correcting code) and let $I$ be an information set for $C$.
A subset $S\subseteq \mathrm{Aut}C$ is a PD-set for $C$ if every $t$-set of coordinate positions can be moved by at least one element of $S$ out of the information set $I$.
The property of having a PD-set for a code is not invariant under isomorphism of codes, it depends on the choice of the information set.

The algorithm of permutation decoding (see \cite{MacS}) uses PD-sets and it is more efficient the smaller the size of a PD-set is.
A lower bound on the size of a PD-set is given in the following theorem and it is due to Gordon \cite{gordon}.

\begin{tm}\label{gordon} 
If $S$ is a PD-set for an $[n,k,d]$ code $C$ that can correct $t$ errors, $r=n-k$, then
$$|S| \geq \left\lceil \frac{n}{r} \left\lceil \frac{n-1}{r-1} \left\lceil \cdots \left\lceil \frac{n-t+1}{r-t+1}  \right\rceil \cdots  \right\rceil  \right\rceil  \right\rceil. $$
\end{tm}

PD-sets for codes do not always exist. Even if they exist, PD-sets are not easy to find, since they depend on the chosen information set of the code. 

Let $A$ be the adjacency matrix of a graph $\Gamma$. Then the full automorphism group of $\Gamma$ is a subgroup of the full automorphism group of the linear code spanned by $A$ over $\mathbb{F}_p$.  
Codes with large automorphism groups are likely to have PD-sets, therefore, we were looking for PD-sets for the codes spanned by adjacency matrices of the DRGs constructed in this paper.

For any of the constructed DRG $\Gamma^{i}_{j}$ from the previous section, let $C^{i}_{j}$ denote the linear code spanned by the adjacency matrix of the graph $\Gamma^{i}_{j}$.
Sizes of the obtained PD-sets (for specific information sets) for some of these codes are given in Table \ref{tb:PD}. For the other codes computation of PD-sets was not feasible.
We denote by $t$ the error correcting capacity of the code, and by $g$ the Gordon bound for the size of the PD-set of a code, from Theorem \ref{gordon}. 
The code $C^{1}_{2}$ is equivalent to the code $C^{2}_{1}$.

\begin{table}[H]
\begin{center} \begin{scriptsize}
\begin{tabular}{|c|c|c|c|c|c|}
\hline
Code $C$ & Parameters $[n,k,d]$ &$Aut(C)$  &$t$ & $g$ &Size of PD-set  \\
\hline
\hline
$C^{1}_{1}$ & [55,10,10] &$S_{11}$ & $4$ &$5$ &$5$\\
$C^{1}_{2}$, $C^{2}_{1}$ & [66,10,20] &$S_{12}$ & $9$ &$15$ &$55$\\
$C^{1}_{5}$ & [330,286,6] &$S_{11}$ & $2$ &$60$ &$420$\\
$C^{1}_{6}$ & [165,120,4] &$S_{11}$ & $1$ &$4$ &$5$\\
$C^{1}_{8}$ & [330,120,8] &$S_{11}$ & $3$ &$7$ &$22$\\
\hline
$C^{3}_{1}$ & [77,20,16] &$M_{22}:Z_2$ & $7$ &$19$ &$110$\\
\hline
$C^{4}_{5}$ & [1771,1540,4] &$S_{23}$ & $1$ &$8$ &$23$\\    
\hline
\hline
\end{tabular} \end{scriptsize}
\caption{\footnotesize PD-sets for codes from constructed DRGs from Mathieu groups}\label{tb:PD}
\end{center}
\end{table}

\subsection{Block designs obtained from a code}

Let $w_i$ denote the number of codewords of weight $i$ in a code $C$ of length $n$. The weight distribution of $C$ is the list $[\left\langle i, w_i: 0 \leq i\leq n\right\rangle]$. The support of a nonzero vector $x=(x_1,...,x_n)\in F_q^n$ is the set of indices of
its nonzero coordinates, i.e. $\mathrm{supp}(x) = \{i | x_i \neq 0\}$. The support design of a code of length $n$ for a given nonzero weight $w$ is the design with points the $n$ 
coordinate indices and blocks the supports of all codewords of weight $w$. 

Here we describe block designs obtained from the code $[176,22,50]_2$ spanned by the adjacency matrix of the graph $\Gamma_2^3$. 

Some remarkable block designs have been constructed from suport designs of codes, for example the 5-designs constructed by Assmus and Mattson in 
\cite{a-m}, and 5-designs constructed by V. Pless in \cite{pless, pless-1}.
In this paper, from the supports of all codewords of the weights of the code $[176,22,50]_2$ we obtain block designs on 176 points on which the finite simple group Higman-Sims acts as the automorphism group.  
The support design for the minimum weight is the very well known Higman-Sims design, denoted by $D_1$ design in Table \ref{tb:des}. 

The weight distribution of the code $[176,22,50]_2$ is:

\noindent
\begin{align*}
&[ \left\langle 0, 1\right\rangle, \left\langle 50, 176\right\rangle, \left\langle 56, 1100\right\rangle, \left\langle 64, 4125\right\rangle,
\left\langle 66, 5600\right\rangle, \left\langle 70, 17600\right\rangle, \left\langle 72, 15400\right\rangle, \left\langle 78, 193600\right\rangle, \\
 &\left\langle 80, 604450\right\rangle, \left\langle 82, 462000\right\rangle, \left\langle 86, 369600\right\rangle, \left\langle 88, 847000\right\rangle, 
\left\langle 90, 369600\right\rangle, \left\langle 94, 462000\right\rangle, \left\langle 96, 604450\right\rangle, \\
 & \left\langle 98, 193600\right\rangle, \left\langle 104, 15400\right\rangle, \left\langle 106, 17600\right\rangle, \left\langle 110, 5600\right\rangle, \left\langle 112, 4125\right\rangle, 
\left\langle 120, 1100\right\rangle, \left\langle 126, 176\right\rangle, \left\langle 176, 1\right\rangle ]
\end{align*}

In Table \ref{tb:des} we describe the results.

\begin{table}[H]
\begin{center} \begin{scriptsize}
\begin{tabular}{|c|c|c||c|c|c|}
\hline
Block design $D$ & Parameters $(v,k,\lambda)$ &$Aut(D)$ & Block design $D$ & Parameters $(v,k,\lambda)$ &$Aut(D)$ \\
\hline
\hline
$D_{1}$ & (176,50,14), $b$=176 & $HS$ &$D_{7}$ & (176,78,37752), $b$=193600 &HS\\
$D_{2}$ & (176,56,110), $b$=1100 &$HS$ &$D_{8}$ & (176,80,124030), $b$=604450 &HS\\
$D_{3}$ & (176,64,540), $b$=4125 &$HS$ &$D_{9}$ & (176,82,99630), $b$=462000 &HS\\
$D_{4}$ & (176,66,780), $b$=5600 &$HS$ &$D_{10}$ & (176,86,87720), $b$=369600 &HS\\
$D_{5}$ & (176,70,2760), $b$=17600 &$HS$ &$D_{11}$ & (176,88,210540), $b$=847000 &HS\\
$D_{6}$ & (176,72,2556), $b$=15400 &$HS$ & &  &\\
\hline
\hline
\end{tabular} \end{scriptsize}
\caption{\footnotesize Block designs from the code $[176,22,50]_2$}\label{tb:des}
\end{center}
\end{table}

\begin{rem}
The action of the group Higman-Sims is transitive on the points of the constructed designs, but in some cases it is not transitive on the blocks of the designs. 
The action that is not transitive on blocks appears in the cases of the designs $D_7$, $D_8$ and $D_{11}$.
\end{rem}

\begin{rem}
In Table \ref{tb:des} we did not include the support designs obtained from the weights 90-176 since they give rise to the complements of the block designs given in Table \ref{tb:des}.
\end{rem}


\vspace*{0.8cm}

\noindent {\bf Acknowledgement} \\
This work has been fully supported by {\rm C}roatian Science Foundation under the project 6732. The authors would like to thank Sven Reichard for pointing out that the graphs $\Gamma_3^1$ and $\Gamma_5^2$ constructed in this paper are not arising from orthogonal arrays.


\end{document}